\documentclass[11pt]{article}
\usepackage[centertags]{amsmath}
\usepackage{graphicx}
\usepackage{amssymb}
\usepackage{amsfonts}
\usepackage{amsthm}
\usepackage{amscd}
\usepackage{fullpage}

\newtheorem{theo}{Theorem}[section]

\newtheorem{coro}[theo]{Corollary}

\newtheorem{defn}[theo]{Definition}

1 \hbadness 9999

\newcommand\Z{{\mathbb Z}}

\newcommand\R{{\mathbb R}}

\begin{document}

\title{A fancy way to obtain the binary digits of $759250125\sqrt{2}$}
\author{Thomas Stoll}
\date{}

\maketitle

%\begin{abstract}
%  R. L. Graham and H. O. Pollak observed that
%  the sequence
%  $$u_1=1,\qquad u_{n+1}=\lfloor \sqrt{2}\,(u_n+\frac{1}{2})\rfloor, \quad n\geq 1,$$
%  has the curious property that the sequence of numbers $(u_{2n+1}-2u_{2n-1})_{n\geq 1}$
%  denotes the binary digits of $\sqrt{2}$. We present an extension
%  of the Graham--Pollak sequence which allows to
%  get -- in a fancy way -- the binary digits of $759250125\sqrt{2}$ and other numbers.
%\end{abstract}

\section{Introduction.}

In the present note we give some ``easily-stated'' recurrences of a special type that generate the binary digits for some ``complicated'' real numbers,
such as the one in the title. The binary digits of any real number $t=(d_1.d_2d_3\cdots)_2$ with $1\leq t<2$ can be calculated by the formula
$$d_n = \lfloor t 2^{n-1}\rfloor -2 \lfloor t2^{n-2}\rfloor,\qquad n\geq 1.$$
We here show a somewhat unexpected, ``fancy'' way to obtain the digits of some special multiples of $\sqrt{2}$, where it is possible to hide
this calculation.

\medskip

This note is structured as follows. In Section~\ref{sec2} we first recall what is known about the Graham--Pollak sequence and its variants, which serve as the motivating examples for the definition of the so-called Graham--Pollak pairs in Section~\ref{sec3}. We then give the general theorem (Theorem~\ref{mtheo1}), which in particular provides a new extension of the original result due to Graham and Pollak. Finally, we give Corollary~3.3 as one surprising example of the general phenomenon.

\section{The Graham--Pollak sequence.}\label{sec2}

As usual, denote by $\lfloor x\rfloor$ the greatest integer less than or equal to $x\in\R$,
and by $\{x\}$  the fractional part of $x$. Define the sequence $(u_n)_{n\geq 1}$ by the recurrence
\begin{equation}\label{GPseq}
  u_1=1,\qquad u_{n+1}=\left\lfloor \sqrt{2}\,\left(u_n+\frac{1}{2}\right)\large\right\rfloor, \quad n\geq 1.
\end{equation}
This sequence, which is also known as the \textit{Graham--Pollak sequence},
 first appeared in a proceedings paper of F. K. Hwang
and S. Lin~\cite{HL69} in the framework of Ford and Johnson's sorting algorithm~\cite{FJ59}.
For the reader interested in the background of the algorithm, an updated exposition
can be found in the third volume of
D.~E. Knuth's \textit{The Art of Computer Programming}~\cite[Ch.~5.3.1, pp.~188]{Kn98}.
The sequence~(\ref{GPseq}) was first investigated from a purely
mathematical point of view by R.~L. Graham and H.~O. Pollak~\cite{GP70}. They found the particularly intriguing fact that
\begin{equation}\label{d}
  d_n=u_{2n+1}-2u_{2n-1}
\end{equation}
gives the $n$th binary digit of $\sqrt{2}=(1.011010100\ldots)_2$.

\medskip

This fact puzzled several authors since then, and it has often
been included as a fun exercise in articles and books mostly
on combinatorial number theory. We mention, for instance,
P. Erd\H os and R.~L. Graham~\cite[p.~96]{EG80},
R. Guy~\cite[Ex.~30]{G88}, R.~L.~Graham, D.~E.~Knuth, and
O.~Patashnik~\cite[Ex.~3.46]{GKP94}. More recent
references are J.-P.~Allouche and J.~Shallit~\cite[Ex.~45, p.~116]{AS03}
and J.~Borwein and D.~Bailey~\cite[p.~62--63]{BB04}. N.~J.~A.~Sloane's
online encyclopedia of integer sequences~\cite{S06} gives eight sequences
which are connected to the Graham--Pollak sequence~(\ref{GPseq}), namely,
A091522, A091523, A091524, A091525, A100671, A100673, A001521, and A004539.

\medskip

Recently~\cite{St05, St06}, the present author found
vast extensions of the Graham--Pollak sequence to parametric
families of recurrences, where the initial value $u_1=1$ is replaced
by $u_1=m$ and the $\sqrt{2}$ in the recurrence is accordingly changed.
However, the sequence is still wrapped in considerable mystery. Indeed, if we do not alter the
$\sqrt{2}$ in the recurrence, but on the other hand, allow only the $1/2$ to vary (if $n$ is odd),
some quite strange things happen: we get the digits of various different multiples
of $\sqrt{2}$, whose digits are seemingly unrelated.
We point out that if we let the $1/2$ vary
for $n$ even instead (cf.~\cite[Theorem~3.3]{St06}), such effects cannot
be observed.

\section{Main Result.}\label{sec3}

In this note we are concerned with the following type of recurrences.

\begin{defn}\label{useq}
{\rm Let $\varepsilon\in\R$ and define the sequence $(v_n)_{n\geq 1}$ by
$$
  v_1=1,\qquad
  v_{n+1}=
 \begin{cases} \lfloor \sqrt{2}\,(v_n+\varepsilon)\rfloor, & \mbox{if }n\mbox{ is odd;} \\
               \lfloor \sqrt{2}\,(v_n+\frac{1}{2})\rfloor, & \mbox{if }n\mbox{ is even.}
 \end{cases}
$$
We call $(\varepsilon,t)$ a {\em Graham--Pollak pair} if the sequence
$$d_n=v_{2n+1}-2 v_{2n-1}, \qquad n\geq 1,$$
represents the binary digits of $t$; that is, $t=(d_1.d_2 d_3\ldots)_2$.}
\end{defn}

\medskip

Note that $(1/2,\sqrt{2})$ is a Graham--Pollak pair according to the original result about the sequence~(\ref{GPseq}).
Our main result is as follows:

\medskip

\begin{theo}\label{mtheo1}
  A list of Graham--Pollak pairs is given by
  $$\{(\varepsilon_i, t_i): \quad 1\leq i\leq 8\},$$
  where
{\footnotesize
  \begin{align*}
     1-\frac{\sqrt{2}}{2}&\leq \varepsilon_1 < \sqrt{2}-1,\quad  &&t_1=\sqrt{2}-1,\\
     \sqrt{2}-1&\leq \varepsilon_2 < \frac{19}{2}\sqrt{2}-13, \quad &&t_2=\frac{11}{8}\sqrt{2}-\frac{5}{8},\\
     \frac{19}{2}\sqrt{2}-13&\leq \varepsilon_3 < \frac{77}{2}\sqrt{2}-54, \quad &&t_3=\frac{45}{32}\sqrt{2}-\frac{19}{32},\\
     \frac{77}{2}\sqrt{2}-54&\leq \varepsilon_4 < \frac{309}{2}\sqrt{2}-218, \quad &&t_4=\frac{181}{128}\sqrt{2}-\frac{75}{128},\\
     \frac{309}{2}\sqrt{2}-218&\leq \varepsilon_5 < \frac{1296121037}{2}\sqrt{2}-916495974, \quad &&t_5=\sqrt{2},\\
     \frac{1296121037}{2}\sqrt{2}-916495974&\leq \varepsilon_6 < \frac{79109}{2}\sqrt{2}-55938, \quad
         &&t_6=\frac{759250125}{536870912}\sqrt{2}-\frac{314491699}{536870912},\\
     \frac{79109}{2}\sqrt{2}-55938&\leq \varepsilon_7 < \frac{5}{2}\sqrt{2}-3, \quad
         &&t_7=\frac{46341}{32768}\sqrt{2}-\frac{19195}{32768},\\
     \frac{5}{2}\sqrt{2}-3&\leq \varepsilon_8 < \frac{\sqrt{2}}{2}, \quad
         &&t_8=\frac{3}{2}\sqrt{2}-\frac{1}{2}.
  \end{align*}
}
\end{theo}

We first comment on a few aspects of the theorem.

\begin{enumerate}
\item[(a)] A surprising feature of Theorem~\ref{mtheo1} is that as $\varepsilon$ varies continuously, the output makes discrete jumps among multiples of $\sqrt{2}$.
Figure~1 illustrates the various intervals for $\varepsilon$ and the corresponding numbers $t$ appearing in Theorem~\ref{mtheo1}.
\item[(b)] The binary digits of $\sqrt{2}$ are obtained for any choice of $\varepsilon$ in the interval
$$[0.4959953\cdots, 0.5012400\cdots).$$
This slightly generalizes the original result of Graham and Pollak with $\varepsilon=1/2$.
\item[(c)] There may well exist Graham--Pollak pairs besides those given in Theorem~\ref{mtheo1}. However, it is easily checked with a computer that the range for admissible values of $\varepsilon$ cannot be too large. For example, for $\varepsilon=0.2928$ we get $d_{3067}=-1$, and for $\varepsilon=0.7073$ we have $d_{2293}=2$. Similar phenomena hold outside these bounds where the values $d_n=-1$ and $d_n=2$ are already obtained for smaller indices $n$. Therefore, possible new pairs can only arise in a very small neighborhood of $\varepsilon=1-\frac{\sqrt{2}}{2}$ or $\varepsilon=\frac{\sqrt{2}}{2}$.
\item[(d)] There is also a connection to normal numbers, which shows that a characterization result for Graham--Pollak pairs is very difficult to obtain. Suppose that there exists $c_1\in\R$ such that
    \begin{equation}\label{normal1}
      \{\sqrt{2}\; 2^{k-1}\}\leq c_1<1,\qquad \mbox{for all } k\geq 1.
    \end{equation}
    Then -- according to~(\ref{conditio}) below -- the interval given for $\varepsilon_1$ can be enlarged to $c_1\left(1-\frac{\sqrt{2}}{2}\right)\leq \varepsilon_1 <\sqrt{2}-1$. Similarly, if there is $c_2\in\R$ such that
    \begin{equation}\label{normal2}
      \{3\sqrt{2}\; 2^{k-2}\}\geq c_2>0,\qquad \mbox{for all } k\geq 1,
    \end{equation}
    then the interval for $\varepsilon_8$ can be enlarged to $\frac{5}{2}\sqrt{2}-3\leq \varepsilon_8<\frac{\sqrt{2}}{2}+c_2\left(1-\frac{\sqrt{2}}{2}\right)$. The inequality~(\ref{normal1}) implies that $\sqrt{2}$ is not normal in base two, and~(\ref{normal2}) implies that $3\sqrt{2}$ is not normal~\cite[Ch.~1.8]{KN74}.
\end{enumerate}

\begin{figure}
\begin{center}
\includegraphics[scale=0.4]{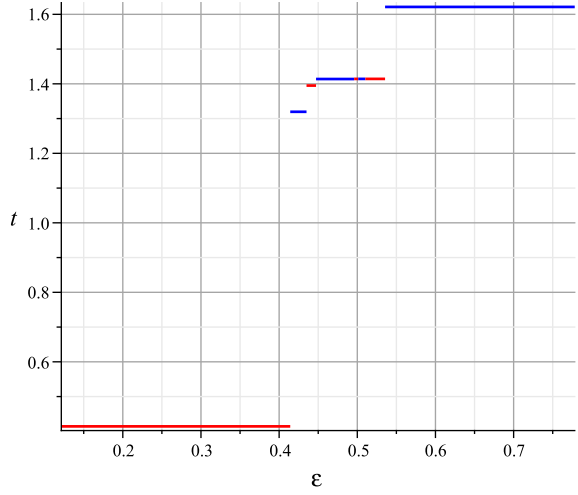}
\caption{The sets $\{(\varepsilon_i,t_i)\}$ for $i=1,\ldots,8$.}
\end{center}
\end{figure}

We conclude with a surprising example, which follows from Theorem~\ref{mtheo1} by the (rather plain) observation that the number $1-\frac{\pi^2}{e^3}=0.5086213\ldots$ lies in the interval given for $\varepsilon_6$.

\begin{coro}\label{cor}
Define the sequence $(w_n)_{n\geq 1}$ by
$$
  w_1=1,\qquad
  w_{n+1}=
 \begin{cases} \lfloor \sqrt{2}\,(w_n+1-\frac{\pi^2}{e^3})\rfloor, & \mbox{if }n\mbox{ is odd;} \\
               \lfloor \sqrt{2}\,(w_n+\frac{1}{2})\rfloor, & \mbox{if }n\mbox{ is even.}
 \end{cases}
$$
Then for $n\geq 31$, $w_{2n+1}-2 w_{2n-1}$ is the $(n+1)$th binary digit of $759250125 \sqrt{2}.$

\end{coro}

\section{Proof of Theorem~\ref{mtheo1}.}

First, let $i\in I:=\{1,2,\ldots,8\}\setminus\{5\}$ and consider the pairs $(\varepsilon_i,t_i)$
in the statement of Theorem~\ref{mtheo1}. Put $$t_i=(\alpha_i \sqrt{2}-\beta_i)\cdot 2^{-l_i}$$
with $\alpha_i,\beta_i,l_i\in\Z$ and $(\alpha_i,2)=1$.
It is easy to verify that $\alpha_i+\beta_i=2^{l_i+1}$ for $i\in I$. Furthermore, let $\xi_{1,i}$ and $\xi_{2,i}$ be the
endpoints of the associated interval for $\varepsilon_i$.
We shall prove that for $\xi_{1,i}\leq \varepsilon_i < \xi_{2,i}$
and $k\geq l_i+2$ we have
\begin{align}
  v_{2k} &= \lfloor t_i 2^{k-2}\rfloor+\gamma_i 2^{k-l_i-2},\label{even}\\
  v_{2k+1} &= \lfloor t_i 2^{k-1}\rfloor+2^k,\label{odd}
\end{align}
where $\gamma_i=2\alpha_i+\beta_i$. This then implies that for $k\geq l_i+3, $$$v_{2k+1}-2 v_{2k-1}=\lfloor t_i 2^{k-1}\rfloor-2 \lfloor t_i 2^{k-2}\rfloor,$$
which is the $k$th binary digit of $t_i$. In the final step we then show that formula~(\ref{odd})
indeed holds true for $0\leq k \leq l_i+1$, which completes the proof.
\medskip

We first use induction to prove that if~(\ref{even}) holds for $k=l_i+2$, then~(\ref{even}) and~(\ref{odd}) hold for $k\geq l_i+2$. Assume the validity of~(\ref{even}). We have to show that
$\left\lfloor \sqrt{2} \left(v_{2k}+\frac{1}{2}\right)\right\rfloor=\lfloor t_i 2^{k-1}\rfloor + 2^k,$
which is equivalent to
$$
  \lfloor t_i 2^{k-1}\rfloor + 2^k  \leq \sqrt{2}\left( \lfloor t_i 2^{k-2}\rfloor+\gamma_i 2^{k-l_i-2}+\frac{1}{2}\right)
<\lfloor t_i 2^{k-1}\rfloor + 2^k+1,
$$
or in other words,
$$  0\leq 2^{k-l_i-1}\left(\beta_i-\frac{\sqrt{2}}{2} \beta_i +\frac{\sqrt{2}}{2}\gamma_i-2^{l_i+1}\right)+\sqrt{2}\; \lfloor \alpha_i \sqrt{2}\; 2^{k-l_i-2}\rfloor
  -\lfloor \alpha_i \sqrt{2} \; 2^{k-l_i-1}\rfloor+\frac{\sqrt{2}}{2}<1.$$
Since $\gamma_i-\beta_i=2\alpha_i$ and $\alpha_i+\beta_i=2^{l_i+1}$ this is the same as
\begin{equation}\label{Ifinal}
  0\leq \{\alpha_i \sqrt{2}\; 2^{k-l_i-1}\}-\sqrt{2} \;\{\alpha_i \sqrt{2}\; 2^{k-l_i-2}\}+\frac{\sqrt{2}}{2}<1.
\end{equation}
Relation~(\ref{Ifinal}) is true since $0\leq \{x\}-\sqrt{2}\;\{x/2\}+\sqrt{2}/2<1$
for all $x\in\R$.

\medskip

Now, assume relation~(\ref{odd}). We have to ensure that
$\left\lfloor \sqrt{2} \left(v_{2k+1}+\varepsilon\right)\right\rfloor=\lfloor t_i 2^{k-1}\rfloor + \gamma_i 2^{k-l_i-1} ,$
or equivalently,
$$
  \lfloor t_i 2^{k-1}\rfloor + \gamma_i 2^{k-l_i-1}  \leq \sqrt{2}\left( \lfloor t_i 2^{k-1}\rfloor+2^k+\varepsilon\right) <\lfloor t_i 2^{k-1}\rfloor + \gamma_i 2^{k-l_i-1}+1.
$$
Here we end up with
\begin{equation}\label{conditio}
  0\leq (1-\sqrt{2}) \{\alpha_i\sqrt{2}\; 2^{k-l_i-1}\}+\sqrt{2}\;\varepsilon <1,
\end{equation}
which is true provided $1-\sqrt{2}/2\leq \varepsilon <\sqrt{2}/2$. This interval includes all of the intervals
$[\xi_{1,i},\xi_{2,i})$ in Theorem~\ref{mtheo1}, and hence there is no additional restriction on $\varepsilon$.

\medskip

It remains to check the initial conditions. This task encompasses some straightforward calculations; we only give the main steps. First, we have to guarantee that~(\ref{even})
is true for $k=l_i+2$. Of course, this crucially depends on the choice of $\varepsilon$. Since $v_n(\varepsilon)$ is non-decreasing for increasing values of $\varepsilon$, there is at most one semi-open real interval $[\bar{\xi}_{1,i},\bar{\xi}_{2,i})$ for $\varepsilon$ such that
\begin{equation}\label{comp}
  v_{2(l_i+2)}=\lfloor t_i 2^{l_i}\rfloor+\gamma_i=\lfloor \alpha_i\sqrt{2}-\beta_i \rfloor+2\alpha_i+\beta_i=\lfloor \alpha_i \sqrt{2} \rfloor +2 \alpha_i.
\end{equation}
We will show that $[\bar{\xi}_{1,i},\bar{\xi}_{2,i})=[\xi_{1,i},\xi_{2,i})$. It is not difficult to crank out a reasonable guess for $\bar{\xi}_{1,i}$ with the help of a computer. In fact,
$v_{2(l_i+2)}$ is a piecewise constant function in $\varepsilon$ with only a finite number of
jump discontinuities. Thus, we can get a close approximation of $\bar{\xi}_{1,i}$ by interval halving.
Furthermore, from Definition~\ref{useq} we see that $\bar{\xi}_{1,i}$ (if it exists) has the form $\frac{c_i}{2}\sqrt{2}-d_i$ for some integers $c_i, d_i\in \Z$. We use \textit{Maple~11} (\texttt{PolynomialTools[MinimalPolynomial]}) to calculate
an approximate minimal polynomial of degree two with ``small'' coefficients to
identify a conjectured value for $\bar{\xi}_{1,i}$. Again, we have to ensure that the value still satisfies~(\ref{comp}).

\medskip

As an illustration, let $i=6$ and consider
$$v_{2(l_i+2)}=v_{62}=v_{62}(\varepsilon),\qquad \lfloor \alpha_i \sqrt{2} \rfloor +2 \alpha_i=2749487923.$$
Figure~2 shows the location of the jumps in the graph of $v_{62}(\varepsilon)$ for $\varepsilon \in [0.40, 0.60]$. By the above procedure we find that $\bar{\xi}_{1,6}$ is ``close'' to $$\xi_{1,6}=1296121037\sqrt{2}/2-916495974=0.5012400\ldots.$$

\begin{figure}
\begin{center}
\includegraphics[scale=0.7]{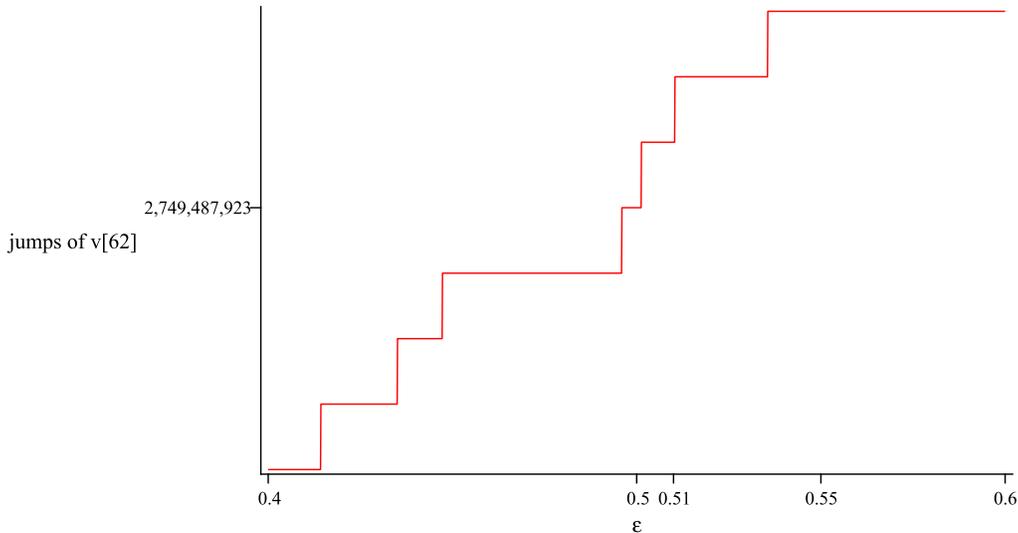}
\caption{The jumps of $v_{62}(\varepsilon)$ for $i=6$ and $0.4\leq \varepsilon \leq 0.6$.}
\end{center}
\end{figure}

\medskip

Once more, we use \textit{Maple} with the \textit{ansatz} $\varepsilon=\xi_{1,i}-\delta$, where $\delta$ denotes a small positive quantity,
to show that $\varepsilon=\xi_{1,i}$ is indeed the smallest value which satisfies~(\ref{comp}). This is a symbolic computation and does not involve high-precision arithmetic.
In a similar fashion, we show that $\bar{\xi}_{2,i}=\xi_{2,i}$. It is important to note that the values of $v_1, v_3, \ldots, v_{2(l_i+1)+1}$ remain unchanged for $\varepsilon\in[\xi_{1,i},\xi_{2,i})$
for every fixed $i\in I$. Moreover, a routine calculation confirms that~(\ref{odd}) is true
for $0\leq k \leq l_i+1$.

\medskip

Finally, we have to treat the case $i=5$, which is less involved than the cases $i\in I$. Here we directly show that
$$v_{2k}=\lfloor t_i 2^{k-2}\rfloor+2^{k-2}\qquad \mbox{and} \quad v_{2k+1} = \lfloor t_i 2^{k-1}\rfloor+2^k$$
for $k\geq 1$, so that we do not have to bother about initial conditions. (We leave the details to the interested reader.)

\medskip

Summing up, we have that the intervals $[\xi_{i,1},\xi_{i,2})$ are disjoint for $i=1,2,\ldots,8$ and completely cover $[1-\frac{\sqrt{2}}{2},\frac{\sqrt{2}}{2})$. This finishes the proof of Theorem~\ref{mtheo1}. \qed

\section*{Acknowledgements.}

The author is a recipient of an APART-fellowship of the Austrian Academy of Sciences at the University of Waterloo, Canada. He wants to express his gratitude to both institutions. The author is also grateful to Jeffrey Shallit for several helpful comments.

\medskip

\noindent \textbf{Thomas Stoll} received his Ph.D. from Graz University of Technology (Austria) in 2004.
Currently, he is a recipient of an APART-fellowship of the Austrian Academy of Sciences at
the University of Waterloo, Canada. When he is not doing mathematics, he is playing the bassoon
and enjoying foreign languages.

\bigskip

\noindent \textit{University of Waterloo, School of Computer Science, Waterloo, ON, Canada, N2L3G1}\\
\textit{tstoll@cs.uwaterloo.ca}

\end{document}